\documentclass[preprint,12pt]{elsarticle}

\usepackage{amssymb}
\usepackage{amsmath}
\usepackage{url}

\DeclareMathOperator*{\argmin}{arg\,min}

\journal{Nuclear Physics B}

\begin{document}

\begin{frontmatter}



\title{Frequency-domain alignment of heterogeneous, multidimensional separations data through complex orthogonal Procrustes analysis} 


\author[ugr]{Michael Sorochan Armstrong} 

\affiliation[ugr]{organization={Computational Data Science Lab, Department of Signal Theory, Telematics and Communications, University of Granada},
            addressline={C/ Periodista Daniel Saucedo Aranda, s/n}, 
            city={Granada},
            postcode={18071}, 
            state={Andalucía},
            country={Spain}}


\begin{abstract}
Multidimensional separations data have the capacity to reveal detailed information about complex biological samples. However, data analysis has been an ongoing challenge in the area since the peaks that represent chemical factors may drift over the course of several analytical runs along the first and second dimension retention times. This makes higher-level analyses of the data difficult, since a 1-1 comparison of samples is seldom possible without sophisticated pre-processing routines. Further complicating the issue is the fact that closely co-eluting components will need to be resolved, typically using some variants of Parallel Factor Analysis (PARAFAC) \cite{armstrong2023parafac2}, Multivariate Curve Resolution (MCR) \cite{parastar2013solving}, or the recently explored Shift-Invariant Multi-linearity \cite{schneide2025unlocking}. These algorithms work with a user-specified number of components, and regions of interest that are then summarized as a peak table that is invariant to shift. However, identifying regions of interest across truly heterogeneous data remains an ongoing issue, for automated deployment of these algorithms \cite{giebelhaus2022untargeted}. This work offers a very simple solution to the alignment problem through a orthogonal Procrustes analysis of the frequency-domain representation of synthetic multidimensional separations data, for peaks that are logarithmically transformed to simulate shift while preserving the underlying topology of the data. Using this very simple method for analysis, two synthetic chromatograms can be compared under close to the worst possible scenarios for alignment.
\end{abstract}

\begin{graphicalabstract}
\end{graphicalabstract}


\end{frontmatter}



\section{Introduction}
\label{sec:intro}
Chromatographic separations separate mixtures of chemicals by exploiting the differences in affinity of each pure component using a functionalised stationary phase. Common technologies are Liquid Chromatography, and Gas Chromatography - both of which may be coupled to a multivariate detector such as a mass spectrometer to enable further identification, and/or deconvolution when the resolving power of the chemical separation step is not sufficient. These coupled instruments, in particular LC-MS and GC-MS have been widely applied to analyze biological samples to interrogate the relative expression of molecules of interest as a function of some known characteristics about the samples themselves. These types of analyses are often referred to as ``un-targeted'' when the goal is to find components that differentiate between interesting population characteristics, as is often the case for poor-understood diseases. When the target molecules are known, the challenges are much less significant - since the points within each chromatogram where they elute can be found by running standards in the absence of the biological matrix.

Through repeated cycles of use, the stationary phase is known to degrade. This leads to the problem of ``peak drift'', and alignment of the signals becomes necessary prior to any kind of multivariate analysis. For multidimensional separations data, this problem becomes two-fold and chemical components that elute within a two-dimensional space may also shift considerably across the two separation modalities.

2D-Correlation Optimized Warping (COW) has previously been proposed as a method to align raw two-dimensional chromatographic signals by Zhang et al. \cite{zhang2008two}, following the work of Tomasi et al. for 1D data \cite{tomasi2004correlation} but relies on 2 parameters: the number of segments and the maximum degree of warping. These parameters must be optimized in order for the analysis to be effective, and for truly heterogeneous data, where two otherwise comparable peaks may exist outside of the pre-defined segments, this technique may not be suitable. This was recently demonstrated by Sorochan Armstrong et al. \cite{sorochan2025alignment} on synthetic data with known characteristics.

This work explores the use of a very simple algorithm to analyze data that has been generated synthetically to simulate the worst possible conditions for peak drift. The images are first transformed using a 2D FFT, and aligned using complex orthogonal Procrustes analysis following vectorization, and down-sampling in the frequency domain. These results suggest its widespread applicability to a variety of different problems involving the alignment of multidimensional separations data, but may require some care when the same compounds are not found in every sample. For this work, we define ``heterogeneous'' as data that is misaligned in a way that may be attributed to different column dimensions or other configuration parameters, but otherwise present a consistent peak order. This would disqualify the use of different stationary phase chemistries that would alter the selectivity of the analysis. 

\subsection{Fast Fourier Transforms}
\label{subsec:fft}

Fourier analyses can be used to encode time, or spatial domain information, such as chromatographic separations data, as a series of complex coefficients that encode the relative amplitude and offset of a series of sinusoidal components through an abstraction of the data that can be analysed directly when there are known experimental characteristics that can be used to factorize the data into orthogonal reconstructions \cite{sorochan2025alignment}. This approach is known as \textit{FFT-ASCA}. Fourier transforms are a linear transformation, that in the one dimensional case can be described as:

\begin{equation}\label{eq:fft}
    x_k = \sum_{n=0}^{N-1} x_n e^{-i \frac{2\pi}{N} kn}\textnormal{, for } k \in \{0...N-1\}
\end{equation}

\noindent for $n \in \{0...N\}$ acquisitions in the time domain where the right hand side of Equation \ref{eq:fft} transforms the data into $k$ frequencies up to $N/2$, or half the total number of acquisitions where $N/2$ is the limiting \textit{Nyquist} frequency \cite{smith1995handbook}. The relationship is linear, and can be inverted to recover the time-domain signal.

\begin{equation}\label{eq:ifft}
    x_n = \frac{1}{N} \sum_{k=0}^{N-1} x_k e^{i \frac{2\pi}{N} kn}\textnormal{, for } n \in \{0...N-1\}.
\end{equation}

The Fast Fourier Transform is an algorithmic implementation of the linear transformation that uses a ``divide and conquer'' strategy to perform the operation in $\mathcal{O}NlogN$ time complexity following the work of Cooley and Tukey \cite{cooley1965algorithm}. In matrix notation, by defining the Discrete Fourier Transform (DFT) matrix \(\mathbf{F} \in \mathbb{C}^{N\times N}\) with entries
\begin{equation}
\mathbf{F} = e^{-i\frac{2\pi}{N} k n}, \quad k,n=0,\dots,N-1,
\end{equation}
we can write the DFT as
\begin{equation}
    x_k = \mathbf{F}x_n,
\end{equation}

\noindent and the inverse as
\begin{equation}
x_n = \frac{1}{N}\mathbf{F}^H\, x_k,
\end{equation}

\noindent where \(\mathbf{F}^H\) denotes the conjugate transpose of \(\mathbf{F}\).

For a two-dimensional signal, consider a matrix 
\(\mathbf{X} \in \mathbb{C}^{M \times N}\) whose elements are \(x_{m,n}\) with 
\(m = 0,1,\dots,M-1\) and \(n = 0,1,\dots,N-1\). The two-dimensional DFT is defined by
\begin{equation}\label{eq:2dfft}
    X_{k,\ell} = \sum_{m=0}^{M-1}\sum_{n=0}^{N-1} x_{m,n}\, e^{-i2\pi \left(\frac{k m}{M} + \frac{\ell n}{N}\right)},
\end{equation}
for \(k = 0,1,\dots,M-1\) and \(\ell = 0,1,\dots,N-1\).

The two-dimensional transform can be computed by applying the one-dimensional transform along each dimension. In matrix notation, if we define the DFT matrices 
\(\mathbf{F}_M \in \mathbb{C}^{M \times M}\) and \(\mathbf{F}_N \in \mathbb{C}^{N \times N}\) with entries
\begin{equation}
[\mathbf{F}_M]_{k,m} = \frac{1}{\sqrt{M}}e^{-i2\pi \frac{km}{M}}, \quad k,m=0,\dots,M-1,
\end{equation}
\begin{equation}
[\mathbf{F}_N]_{l,n} = \frac{1}{\sqrt{N}}e^{-i2\pi \frac{ln}{N}}, \quad l,n=0,\dots,N-1,
\end{equation}
then the two-dimensional DFT of \(\mathbf{X}\) is given by
\begin{equation}
\mathbf{X}_f = \mathbf{F}_M\, \mathbf{X}\, \mathbf{F}_N^H,
\end{equation}
where the superscript \(H\) denotes the conjugate transpose. (In the special case where \(M=N\), we often write this as \(\mathbf{X}_f = \mathbf{F}\,\mathbf{X}\,\mathbf{F}^H\).)

The inverse two-dimensional Fourier transform is given by
\begin{equation}
\mathbf{X}_{m,n} = \frac{1}{MN}\sum_{k=0}^{M-1}\sum_{\ell=0}^{N-1} X_{k,\ell}\, e^{i2\pi\left(\frac{km}{M}+\frac{\ell n}{N}\right)},
\end{equation}
\noindent or in matrix form:
\begin{equation}\label{eq:ifft2}
\mathbf{X}_{m,n} = \mathbf{F}_M^H\, \mathbf{X}_f\, \mathbf{F}_N.
\end{equation}

\subsection{Orthogonal Procrustes Analyses}

Orthogonal Procrustes analysis is a method that maps matrix $\mathbf{A}$ onto matrix $\mathbf{B}$ by finding the optimal orthogonal subspace spanned by the matrix $\mathbf{\Omega}$ within a common number of dimensions. The problem is formalized as follows:

\begin{equation}
    \argmin_{\mathbf{\Omega}}||\mathbf{\Omega} \mathbf{A} - \mathbf{B}||_F^2
\end{equation}

\noindent and the optimal, closed-form solution exploits the decomposition of the matrix $\mathbf{M} = \mathbf{B}\mathbf{A}^T$ using the left and right singular vectors $\mathbf{U}$ and $\mathbf{V}$ \cite{hurley1962procrustes}:

\begin{equation}
    \mathbf{\Omega} = \mathbf{U}\mathbf{V}^T.
\end{equation}

For \textit{generalized} Procrustes analysis, variance information encoded in the singular values themselves can be incorporated, to stretch or compress particular $B$ to map to $A$. For the sake of preserving quantitative information, only orthogonal Procrustes analysis is used in this work. Further details and applications of Procrustes analysis can be found in the work by Kucheryavskiy et al \cite{kucheryavskiy2023procrustes}.

\subsection{Logarithmic distortion}

We generate a series of blobs randomly in an image with width $W$, and height $H$ and distort the coordinates of these blobs logarithmically to preserve the underlying topology of the data, while ensuring continuity (i.e. the coordinates of each blob must not exceed the width or height of the image). The coordinates are defined:

\begin{equation}
    x \in [1, W] \quad \text{and} \quad y \in [1, H].
\end{equation}

To generate the modified logarithmic mapping as a function of the original coordinates, we introduce a blending factor \(\alpha\) (with \(0 \le \alpha \le 1\)). The final transformed coordinates \(x'\) and \(y'\) which are then defined as
\begin{equation}\label{eq:log1}
    x' = (1-\alpha)x + \alpha\left(1 + (W-1)\frac{\ln(x)}{\ln(W)}\right),
\end{equation}

\begin{equation}\label{eq:log2}
y' = (1-\alpha)y + \alpha\left(1 + (H-1)\frac{\ln(y)}{\ln(H)}\right).
\end{equation}

This distortion insures that there is no ``peak inversion'', where the coordinates of the blobs in the distorted image would map inversely to their original coordinates. In the absence of mass spectral, or other identifying information - the fact that the data is restricted so that it cannot pass through itself, is seen as a reasonable limitation, especially in the context of its underlying topology. 

\subsection{FFT-Orthogonal Procrustes Analysis}

A distorted image, representing multidimensional separations data is defined as $\mathbf{X}_d$ with width and height $W$ and $H$, respectively. The target image is defined as $\mathbf{X}_t$, which $\mathbf{X}_d$ will map to. The 2D Fourier transforms of the images are first taken, and the multidimensional coefficients are vectorized for subsequent analysis for $\tilde{\mathbf{X}}_d = \mathbf{F}\mathbf{X}_d\mathbf{F}^H$ and $\tilde{\mathbf{X}}_t = \mathbf{F}\mathbf{X}_t\mathbf{F}^H$. The vectorized forms are defined as $x_d \in 1, W \times H$ and $x_t \in 1, W\times H$ following $\tilde{x}_d = vec(\tilde{\mathbf{X}}_d)$ and $\tilde{x}_t = vec(\tilde{\mathbf{X}}_t)$. Following this, the complex coefficients are calculated using the uni-dimensional FFT, denoted by the tilde: $\tilde{x}_d$, and $\tilde{x}_t$. The orthogonal Procrustes problem is then denoted as:

\begin{equation}
    \argmin_\mathbf{\Omega} ||\mathbf{\Omega} \tilde{x}_d - \tilde{x}_t||^2
\end{equation}

\noindent and the complex-valued $\Omega$ generalizes easily from the real-valued solution to the orthogonal Procrustes problem, from the decomposition of the other product of $\tilde{x}_d$ and $\tilde{x}_t$, through the \textit{Hermitian}, or complex transpose:

\begin{equation}
    \tilde{\mathbf{M}} = \tilde{x}_t \tilde{x}_d^H.
\end{equation}

Similarly, given that $\tilde{\mathbf{M}}$ can be decomposed following $\tilde{\mathbf{M}} = \mathbf{U}\mathbf{\Sigma}\mathbf{V}^H$, the optimal value for $\mathbf{\Omega}$ is calculated using the left and right singular vectors:

\begin{equation}
    \mathbf{\Omega} = \mathbf{U}\mathbf{V}^H.
\end{equation}

The ``distorted'', frequency-domain representation of $\tilde{x}_d$ is mapped to the target $\tilde{x}_t$ following:

\begin{equation}
    \hat{x}_t = \mathbf{\Omega}\tilde{x}_d
\end{equation}

\noindent and then $\hat{x}_t$ is folded back into matrix form as: $\hat{\mathbf{X}}_t$ before being transformed back into the spatial domain following Equation \ref{eq:ifft2}. We assume that $\mathbf{\Omega}$ is unitary such that $\mathbf{\Omega}\mathbf{\Omega}^H = \mathbf{I}$. This extension is trivial and can be proven using similar mathematics to the derivation of the original orthogonal Procrustes problem.

\section{Materials and Methods}

All data were synthetically generated in MATLAB R2024a, and analysed locally in the same environment. The calculations were performed on system running Ubuntu 22.04.4 LTS with an Intel i9-1400K 32 core CPU, two parallel NVIDIA GeForce RTX 4070 GPUs and 128 GB of RAM. All code used in this analysis is available online at: \url{https://github.com/mdarmstr/fft_procrustes_demo}.

\section{Results and Discussion}

\begin{figure}
    \centering
    \includegraphics[width=0.95\linewidth]{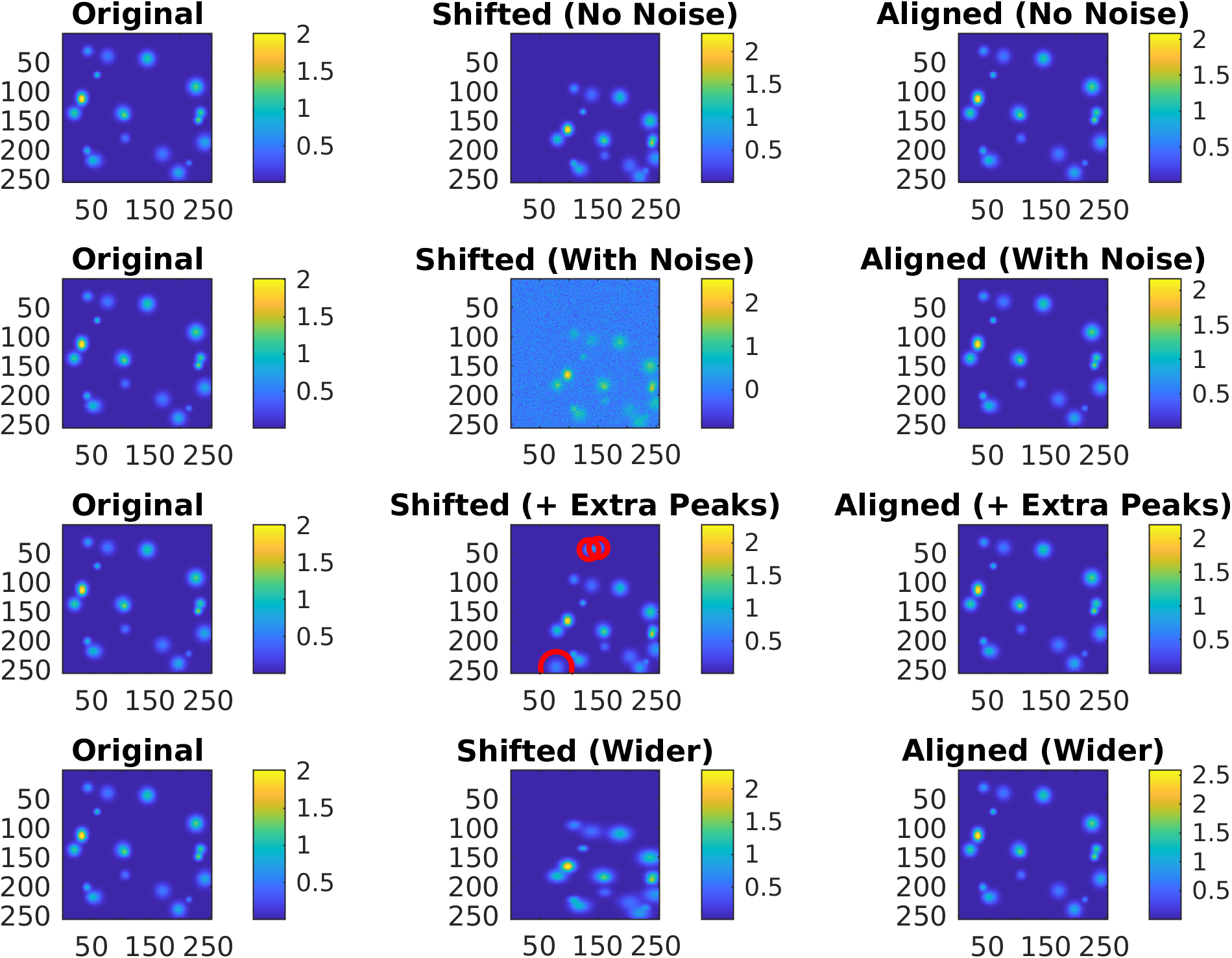}
    \caption{Performance of the FFT-Procrustes analysis. Row 1 shows the original, shifted, and aligned image without any noise. Row 2 shows the original, shifted and aligned image with 10\% Gaussian noise as a function of the maximum peak intensity. Row 3 shows the results with extra peaks added to the shifted image, and on row 4 the shifted peaks have their widths increased by 50\%.}
    \label{fig:analysis}
\end{figure}

4 test cases were examined, to evaluate the performance of the algorithm under the worst possible situations for peak drift. Each following the logarithmic transformations according to Equations \ref{eq:log1} and \ref{eq:log2}. The results are shown in Figure \ref{fig:analysis}. All blobs were generated according to a 2-dimensional Gaussian distribution, but were not restricted so as to not overlap in either the original or distorted images. The dimensions of each image were set as $H, W = 256$ and 20 blobs were generated. The number of frequency components was down-sampled by a factor of 4, to facilitate the analysis directly on the device's volatile storage. An $\alpha$ of 0.5 was used to distort the position of each Gaussian blob according to Equations \ref{eq:log1} and \ref{eq:log2}.

The algorithm performs well regardless of a significant amount of noise, the presence of 3 extra peaks, and even when the peak shapes themselves are not consistent with the original image. Within both the original and distorted images there is a considerable amount of peak overlap. It is worth noting however, that in the instance where the extra peaks are introduced to the shifted image they do not appear in the aligned image. For laboratory practitioners, this drawback may be avoided by the use of pooled samples - wherein each chromatogram is aligned against an image containing all possible chemical components.

For a quantitative measure of performance, the cosine correlation coefficient between the original as shown below:
\begin{equation}
  \rho_{\cos}(\mathbf{x}, \mathbf{y}) = \frac{\mathbf{x} \cdot \mathbf{y}}{\|\mathbf{x}\|\,\|\mathbf{y}\|} 
= \frac{\sum_{i=1}^{n} x_i y_i}{\sqrt{\sum_{i=1}^{n} x_i^2}\sqrt{\sum_{i=1}^{n} y_i^2}}.  
\end{equation}

\noindent where $i \in n$ describes each entry in the unfolded vector. The results comparing the aligned datasets with the originals were performed and summarized in Table \ref{tab:cos}:
\begin{table}[ht]\label{tab:cos}
\begin{center}

\begin{tabular}{l c}
\hline
\bf Variant & \bf Cosine Correlation \\
\hline
Original (self) & 1.0000 \\
No Noise & 0.9999 \\
With Noise & 0.9999 \\
+ Extra Peaks & 0.9999 \\
Wider & 0.9999 \\
\hline
\end{tabular}
\end{center}
\caption{Tabular summary of cosine correlation coefficients between the aligned and original data, demonstrating good agreement overall.}
\end{table}

\section{Conclusions}

This method introduces a very simple pre-processing strategy that can be used to analyze multidimensional separations data on a pixel-wise basis. It may also further enable the use of simple methods for deconvolution, such as PARAFAC without the need for additional model constraints such as those used in PARAFAC2 for the analysis of GC-MS data \cite{kiers1999parafac2}. Furthermore, this method can also be used in the absence of a multivariate detector such as a mass spectrometer, as this has been demonstrated using only image data that would correspond to an analysis performed using a univariate detector.

This method highlights that while some phenomena may be non-linear in the time or spatial domains, a simple linear solution may be present through an analysis in the frequency domain.

A drawback of the technique is the requirement of a large amount of memory to calculate the rotation matrix $\mathbf{\Omega}$. This may be avoided by caching additional memory onto non-volatile storage devices, but will slow down the analysis considerably. 

\section{Acknowledgements}

Michael Sorochan Armstrong has received funding from the European Union's Horizon Europe Research and Innovation Program under the Marie Skłodowska- Curie grant agreement Meta Analyses of Heterogeneous Omics Data (MAHOD) no. 101106986. This work was supported by grant no. PID2023-1523010B-IOO (MuSTARD), funded by the Agencia Estatal de Investigación in Spain, call no. MICIU/AEI/10.13039/501100011033, and by the European Regional Development Fund, with MICIU, that comes from Ministerio de Ciencia, Innovación y Universidades.

The author would like to thank Drs. Edoardo Saccenti and Jos\'e Camacho from the Universities of Wageningen and Granada for their feedback and encouragement, and Dr. Michelle Corbally from Los Alamos National Laboratories, USA for her thoughtful commentary at the 16th Multidimensional Chromatography Workshop in Li\`ege, Belgium that was the inspiration for this work.

\bibliographystyle{ieeetr}
\bibliography{cas-refs}

\end{document}